\documentclass[12pt]{article}
\usepackage{amsmath,amssymb,amsfonts}

\newcommand{\D}{\displaystyle}
\newcommand{\R}{\mathbb{{R}}}

\newcommand{\N}{\mathbb{{N}}}

\newcommand{\lo}{\longrightarrow}

\newcommand{\ov}{\overline}
\newcommand{\st}{\stackrel}

\newcommand{\para}{\paragraph}
\begin{document}
\begin{center}
{\bf On Conley's Fundamental Theorem of Dynamical Systems}\\
[1cm] M.R. Razvan \\
{\it\small (To the memory of C. Conley)}
\end{center}

\vspace*{5mm}
\begin{abstract}
In this paper, we generalize Conley's fundamental theorem of
dynamical systems in Conley index theory.  We also conclude the
existence of regular index filtration for every Morse
decomposition. 

\end{abstract}

\section{Introduction}
Charles Conley is mostly known for his fundamental theorem of
dynamical systems and his homotopy index theory \cite{C}. In the
latter, he proved that every continuous flow on a compact metric
space admits a Lyapunov function which is strictly decreasing
along the orbits off the chain recurrent set. This result has been
developed by Franks for homeomorphisms \cite{F} and Hurley for
noncompact metric spaces \cite{H1}-\cite{H4}. In the former,
Conley defined a homotopy invariant for any isolated invariant set
for a continuous flow. This invariant gives some invaluable
information about the behavior of the isolated invariant set.
This paper concerns a combination of these two masterpieces.
Indeed we generalize Conley's fundamental theorem of dynamical
systems in Conley index theory. We also conclude the existence of
regular index filtration for every Morse decomposition. The
generalized Morse inequalities is a consequence of the existence
of regular index filtration.

\section{Conley index theory}
We start this section with some basic definitions. Let $\varphi^t$
be a continuous flow on a metric space $X$. An isolated invariant
set is a subset $S \subset X$ which is the maximal invariant set
in a {\bf compact} neighborhood of itself. Such a neighborhood is
called an isolating neighborhood. A Morse decomposition for $S$ is
a collection $\{M_i\}^n_{i=1}$ where each $M_i$ is an isolated
invariant subset of $S$ and for all $x\in S-\bigcup_{i=1}^n M_i$
there exist $i,j\in \{1,\cdots , n\}$ such that $i>j$, $\alpha
(x) \in M_i$ and $\omega (x)\in M_j$. A pair $(A,A^*)$ of subsets
of $S$ is called an attractor-repeller pair if $\{A,A^*\}$ is a
Morse decomposition for $S$ i.e. $\alpha (x) \in A^*$ and $\omega
(x)\in A$ for every $x\in S-(A\cup A^*)$.

Let $S$ be an isolated invariant set with an isolating
neighborhood $V$ and a Morse decomposition $\{M_i\}^n_{i=1}$. In
\cite{S}, it is proved that if $\varphi^{[0,+\infty)} (x) \subset
V$, then $\omega (x)\subset M_i$ for some $1\leq i\leq n$.
Similarly if $\varphi^{(-\infty ,0]} (x) \subset V$, then $\alpha
(x)\subset M_i$ for some $1\leq i\leq n$. Now for $j=0,\cdots ,
n$, we define
\begin{center}
$I^+_j=I^+_j(V)=\{x\in V | \varphi^{[0,\infty)} (x) \subset V ,\
\omega (x)\subset M_{j+1}\cup \cdots \cup M_n\}$\\
$I^-_j=I^-_j(V)=\{x\in V | \varphi^{(-\infty ,0]} (x) \subset V,
\  \alpha (x)\subset M_1\cup \cdots \cup M_j\}$ \\
$S_j^*=\{x\in S|\omega (x)\subset M_{j+1}\cup \cdots \cup M_n\}\
, \ S_j=\{ x\in S | \alpha (x)\subset M_1\cup \cdots \cup M_j\}$. \\
\end{center}
 Moreover if
$(A,A^*)$ is an attractor-repeller pair for $S$, we set

$I^+_{A^*}=\{x\in V | \varphi^{[0,\infty)} (x) \subset V ,\
\omega (x)\subset A^*\}\ ,\ I^-_A=\{x\in V | \varphi^{(-\infty
,0]} (x) \subset V ,\  \alpha (x)\subset A\}$

In \cite{S} and \cite{CZ}, it is proved that  $I_j^+$ and $I_j^-$
are compact and $(I_j^+, I_j^-)$ is an attractor-repeller pair for
$S$. This fact allows us to prove our results for an
attractor-repeller pair and then extend them to every Morse
decomposition.

In order to define the concept of index pair , we follow \cite{RS}
and \cite{S}.
 Given a compact pair $(N,L)$ with $L\subset N\subset X$, we define the induced semi-flow
on $N/L$ by $$\varphi_{\sharp}^t: N/L \lo N/L ,  \ \
\varphi_{\sharp}^t (x)=\left\{
\begin{array}{ll}
\varphi^t (x) & \text{if} \ \varphi^{[0,t]} (x) \subset N-L\\
\left[ L \right] & \text{otherwise.}
\end{array} \right. $$
In \cite{RS} it is proved that $\varphi^t_{\sharp}$ is continuous
if and only if

i) $L$ is positively invariant relative to $N$, i.e. $$x\in L, \
t\geq 0,\  \varphi^{[0,t]} (x) \subset N \Rightarrow
\varphi^{[0,t]} (x) \subset L.$$

ii) Every orbit which exits $N$ goes through $L$ first, i.e.
$$x\in N,\varphi^{[0,\infty)} (x) \not\subset N \Rightarrow
\exists_{t\geq 0} \ \text{with} \ \varphi^{[0,t]} (x) \subset N,
\varphi^t (x) \in L,$$ or equivalently if $x\in N-L$ then there is
a $t>0$ such that $\varphi^{[0,t]} (x) \subset N$.

\para{Definition.} An index pair for an isolated invariant set
$S\subset X$ is a compact pair $(N,L)$ in $X$ such that
$\overline{N-L}$ is an isolating neighborhood for $S$ and the
semi-flow $\varphi^t_{\sharp}$ induced by $\varphi^t$ is
continuous.

 In \cite{C}, \cite{CZ}, \cite{RS}
 and \cite{S} it has been shown that every  isolated invariant set $I$ admits
an index pair $(N,L)$ and the homotopy type of the pointed space
$N/L$ is independent of the choice of the index pair. The Conley
index of $S$ is the homotopy type of $(N/L,[L])$.\\
{\bf Note.} We shall not distinguish between $N-L$ and
$N/L-\{[L]\}$.

\para{Definition.} An index pair $(N,L)$ is called regular if the
exit time map defined by  $$\tau_+:N \lo [0,+\infty) , \ \
\tau_+(x)=\left\{
\begin{array}{ll}
sup\{t|\varphi^{[0,t]} (x) \subset N-L\} & \text{if} \ x\in N-L,\\
0 & \text{if}\ x\in L,
\end{array} \right. $$
is continuous. For every regular index pair $(N,L)$, we define the
induced semi-flow on $N$ by $$ \varphi ^t _{\natural} :N\times \R
^+ \lo N,\ \ \varphi^t _{\natural}(x)=\varphi ^{min\{t,\tau_+
(x)\}} (x)$$.

\para{Proposition 2.1.} If $(N,L)$ be a regular index pair for a
continuous flow $\varphi ^t$, then $L$ is a neighborhood
deformation retract in $N$. In particular, the natural map $\pi :
N \lo N/L$ induces an isomorphism $H_*(N,L)\cong H_*(N/L,[L])$. \\
{\bf Proof.} Consider the induced semi-flow $\varphi_{\natural}$
on $N$ and the neighborhood $U:=\tau _+ ^{-1} [0,1]$ of $L$. Now
$\varphi_{\natural}|_{U\times [0,1]}$ gives the desired
deformation retraction. $\square$

It is not true that $H_*(N,L)\cong H_*(N/L,[L])$ for every index
pair $(N,L)$. For this reason, it is common to use the
Alexander-Spanier cohomology for which the above isomorphism is
always valid \cite{C}, \cite{CZ}. However it is easier to deal
with more familiar homology and cohomology theories. Therefore we
are interested in regular index pairs. the following  result
provides a criterion for regularity of an index pair. The reader
is referred to \cite{S} for the details about regular index pairs
and the proof of the following useful criterion.

\para{Proposition 2.2.} An index pair $(N,L)$ is regular provided
that $\varphi ^{[0,t]}(x)\not \subset \ov{N-L}$ for every $x\in L$
and $t>0$.

\para{Definition.} Let $S$ be an isolated invariant set with a
Morse decomposition $\{M_i\}_{i=1}^n$. An index filtration is a
sequence $N_0\subset N_1\subset \cdots N_n$ of closed subsets of
$X$ such that $(N_k,N_{k-1})$ is an index pair for $M_k$ for
every $1\leq k\leq n$. When each $(N_k,N_{k-1})$ is regular, then
the above filtration is called a regular index filtration.

It is not hard to see that if $N_0\subset N_1\subset \cdots N_n$
is an index filtration, then  $(N_n, N_0)$ is an index pair for
$S$. Moreover if the filtration is regular, then $(N_n, N_0)$ is
a regular index pair for $S$. It is well-known that every Morse
decomposition admits an index filtration and the generalized Morse
inequalities is a consequence of this result \cite{CZ}. In the
next section, we show that every Morse decomposition admits a
regular index filtration.

\section{Conley's Fundamental Theorem}
In this section we generalize Conley's Fundamental Theorem of
Dynamical Systems in Conley index theory. The following lemma
will play a crucial role in the continuity of Conley's Lyapunov
function.

\para{Lemma 3.1.} Let $S$ be an isolated invariant set with an
attractor-repeller pair $(A,A^*)$, an index pair $(N,L)$ and the
isolating neighborhood $V=\ov{N-L}$. If $B$ is a compact subset of
$N/L -I^+ _{A^*} $ and $U$ is a neighborhood of $[L]\cup I^- _A
$, then there exists $T\in \R ^+$ such that $\varphi _{\sharp }
^{[T,+\infty)} (B) \subset U$.

\para{Proof.} We may assume that $U$ is a compact neighborhood of
$[L]\cup I^- _A $ with $U \cap I^- _A =\varnothing$. Now suppose
that there are $x_n \in B$ and $t_n \lo \infty$ such that
$\varphi _{\sharp}^{t_n}(x_n)\not \in U$. Since $B$ is compact,
we may choose $x_n$'s so that $x_n \lo x \in B$. It is easy to
see that $\varphi ^{[0,+\infty)}(x)\subset N-L$. Since $B\cap
I^+_{A^*}=\varnothing$, we have $\omega (x)\subset A$, hence
there is a $t\in \R ^+$ such that $\varphi ^{[t,+\infty)} (x)\in
\st{\circ}{U} $.  Since $x_n \lo x$, there are $t_n \in \R ^+$
such that $t_n \lo +\infty$,  $ \varphi
^{[t,t_n]}_{\sharp}(x_n)\subset \st{\circ}{U}$ and $t_n >t $ for
every $n \in \N$. It follows that there exists $t'_n\in [t,t_n]$
such that $\varphi ^{[t,t'_n]}\subset U,\ \varphi ^{t'_n}(x_n)\in
\partial U$ and $t'_n -t\lo +\infty$. (The latter is true since
 $\varphi ^{[t,+\infty)}(x_n)\subset \st{\circ}{U}$ and $x_n \lo
x$.) Therefore the sequence $\varphi ^{t'_n}(x_n)$ has a limit
point $y \in \partial U$ such that $\varphi ^{(-\infty ,
0]}(y)\subset U\cap (N-L)$ and $y\in \omega (B)\subset S$. Thus
$\alpha (y) \subset A$ which follows that $y\in A$. This
contradicts $y\in
\partial U$. $\square$

\para{Theorem 3.2.} Let $S$ be an isolated invariant set with an
attractor-repeller pair $(A,A^*)$ and an index pair $(N,L)$.
 There exists a continuous function $g:N/L \lo [0,1]$ such that\\
(i) $g^{-1} (0)=[L]\cup I^- _A$ and $g^{-1} (1)= I^+ _{A^*}$.\\
(ii) $g(\varphi ^t (x))<g(x)$ for every  $x \not \in [L]\cup  I^-
_A \cup I^+ _{A^*}$ and $t \in \R ^+$.

\para{Proof.} Let $\rho : N/L \lo [0,1]$ be a continuous function
with $\rho ^{-1}(0)=[L]\cup I^-_A$ and $\rho ^{-1}(1)=I^+_{A^*}$.
We define $h:N/L \lo [0,1]$ by $h(x)=\D{\sup_{t\geq 0}}\
\rho(\varphi ^t_{\sharp}(x))$. It is not hard to see that
$h^{-1}(0)=[L]\cup I^-_A,\ h^{-1}(1)=I^+_{A^*}$. We show that $h$
is upper semi-continuous. For every $x\in N/L$ and $\epsilon
>0$ there is a $t\in \R ^+$ such that $\rho (\varphi ^t
_{\sharp}(x))>h(x)-\epsilon $. Now there is a neighborhood $U$ of
$x$ such that $\varphi ^t _{\sharp}(y)>h(x)-\epsilon$ for every
$y\in U$. Therefore $h(y)>h(x)-\epsilon$ for every $y\in U$ which
proves the upper semi-continuity of $h$. As a result, $h$ is
continuous in $h^{-1}(1)$. Now suppose that $x\not \in
h^{-1}(1)=I^+_{A^*}$ and $\epsilon <1-h(x)$. If we set $B=\rho
^{-1}[0,h(x)+\epsilon ]$ and $U=\rho ^{-1}[0,h(x)+\epsilon )$ in
the above lemma, we obtain a $T\in \R ^+$ with $\rho (\varphi
^t_{\sharp}(y))<h(x)+\epsilon $ for every $y\in B$ and $t\geq T$.
Now by continuity of $\varphi _{\sharp}$, there exists an open
set $V\subset N/L$ such that $x\in V$ and $\rho (\varphi
^t_{\sharp}(y))<h(x)+\epsilon $ for every $t\in [0,T]$ and $y\in
V$, Therefore  $h(y)<h(x)-\epsilon $ for every $t\in [0,T]$ and
$y\in U\cap V$ which shows that $h$ is lower semi-continuous in
$x$. Now it is not hand to check that $g:=\int_0^{+\infty} e^{-t}
f(\varphi^t_{\sharp} (x))dt$ is the desired function \cite{S}.
$\square$

\para{Theorem 3.3.} Let $S$ be an isolated invariant set with an
index pair $(N,L)$ and a Morse decomposition $\{M_i\}_{i=1}^n$.
There is  a continuous function $g:N/L \lo [0,N+1]$ such
that\\ (i) $g^{-1}(0)=[L]$ and $g(M_i)=i$ for every $1\leq i\leq n$.\\
(ii) If $x \in N-L-\bigcup _{i=1}^n M_i$ and $t>0$, then
$g(\varphi ^t_{\sharp}(x))<g(x)$.\\
{\bf Proof.} Consider the attractor-repeller pairs $(I_j, I_j^*)$
for $0\leq j\leq n$. By Theorem 3.2, there are continuous
functions $g_i:N/L\lo [0,1]$ with  $g_i^{-1}(0)=[L]\cup I_j^-,\
 g_i^{-1}(1)=I_j^+$ and
 $g_j(\varphi^t(x))<g_j(x)$ for every
 $x\not \in [L]\cup I_j^-\cup I_j^+$. Now $g:=g_0+\cdots +g_n$ is the
 desired function. $\square$

\para{Corollary 3.4.} Every Morse decomposition admits a regular
index filtration. \\
{\bf Proof.} Let $f$ be the above Lyapunov function. If we set
$N_k:=\pi ^{-1}(f^{-1}[0,k+1/2])$ for $0\leq k\leq n$ and
$N_n:=N$, then by Proposition 2.2, $(N_k, N_{k-1})$ is a regular
index pair for $M_k,$ for every  $1\leq k\leq n$. $\square$

\para{Definition.} Let $\varphi ^t$ be a continuous flow on a
compact metric space $X$. An $\epsilon$-chain for $\varphi ^t$ is
a sequence $x_0,\cdots ,x_n$ in $X$ and $t_1,\cdots ,t_n$ in
$[1,+\infty)$ such that $d(\varphi ^{t_i}(x_{i-1}))<\epsilon $. A
point $x\in X$ is called chain recurrent if for every $\epsilon
>0$, there is an $\epsilon $-chain $x=x_0,\cdots ,x_n =x$. The
set of all chain recurrent points for $\varphi ^t$ is denoted by
$R(\varphi ^t)$.

It is not hard to Check that $R(\varphi ^t)$ is a closed invariant
subset of $X$ containing the non-wandering set $\Omega (\varphi
^t)$. In \cite{C} and \cite{Ro}, it has been shown that $R(\varphi
^t|_{R(\varphi ^t)})=R(\varphi ^t)$ and $R(\varphi ^t)=\bigcap
(A\cup A^*)$ where the intersection is taken over all
attractor-repeller pairs $(A,A^*)$ in $X$. It is also known that
the number of all attractor-repeller pairs in a compact metric
space is at most countable.

\para{Theorem 3.5.} Let $S$ be an isolated invariant set with an
index pair $(N,L)$. Then there is a continuous function $g:N/L\lo
[0,1]$ such that\\
(i) $g^{-1}(0)=[L]$ and $g(\varphi ^t(x))\leq g(x)$ for every
$x\in
N/L$ and $t\geq 0$.\\
(ii)if $x \in N-L-R(\varphi ^t |_I)$ and $t\geq$, then $g(\varphi
^t_ {\sharp}(x))<g(x)$.\\
{\bf Proof.} Let $\{(A_i,A_i^*)\}_{i=1}^{\infty}$ be the sequence
of all attractor-repeller pairs in $S$ including $(\varnothing,
S)$ and $(S,\varnothing)$. Now by the Theorem 3.2, there are
continuous functions $g_i:N/L\lo [0,1]$ such that
$g_i^{-1}(0)=[L]\cup I^-_{A_i},\ g_i^{-1}=I^+_{A^*_i}$ and
$g_i(\varphi ^t(x))<g_i(x)$ for every $t\in \R ^+$ and $x\not \in
[L]\cup I^+_{A^*_i}\cup I^-_{A_i}$. Now $g=\sum _{i=1}^{\infty}
2^{-i} g_i$ is the desired function. $\square$

The above result can be considered as a  generalization of
Conley's Fundamental Theorem of Dynamical Systems. A similar
result for maps and semi-flows can be obtained by following
\cite{FR}, \cite{H3}, \cite{H4} and \cite{Sz}.

\para{Acknowledgment.} The author would like to thank Institute for
studies in theoretical Physics and Mathematics, IPM ,  for
supporting this research.

\noindent  Institute for Studies in Theoretical Physics and
Mathematics \\ P.O.Box: $19395-5746$, Tehran, IRAN. \\
e-mail: razvan@karun.ipm.ac.ir , Fax: 009821-2290648
\end{document}